\documentclass[11pt]{article}
\usepackage{amsmath,amssymb,theorem}
\usepackage[hmargin=3cm,vmargin={3.5cm,4cm}]{geometry}
\newtheorem{theorem}{Theorem}[section]
\newtheorem{remark}{Remark}[section]
\newtheorem{lemma}[theorem]{Lemma}
\newtheorem{proposition}{Proposition}[section]
\numberwithin{equation}{section}
\numberwithin{theorem}{section}

\newcommand{\R}{\mathbb R}

\begin{document}

\title{Confinement of a hot temperature patch in the modified SQG model.}

\author{
Roberto Garra \footnote{Email address:
roberto.garra@uniroma1.it},	\\
    Department of Statistical Sciences, Sapienza University of Rome,\\
       Piazzale Aldo Moro 5, 00185 Rome, Italy.\\
}

\date{}
\maketitle

\begin{abstract}
In this paper we study the time evolution of a temperature patch in
$\mathbb{R}^2$ according to the modified Surface
Quasi-Geostrophic Equation (SQG) patch equation. In particular we give a temporal
estimate on the growth of the support, providing a rigorous proof 
of the confinement of a hot patch of temperature in absence of external 
forcing, under the quasi-geostrophic approximation.
\end{abstract}
\

\bigskip\noindent

\section{Introduction}
\label{sec:1}

In the present paper we study the time evolution of a patch of
temperature in $\mathbb{R}^2$ according to the modified surface quasi-geostrophic (SQG) equation, considered in recent years by 
different authors (see e.g. \cite{kisel} and the references therein). We recall the basic equations of the $\alpha$-patch model: let
$\theta (x,t)$, $x \in \R^2$ be the solution of the equation
\begin{equation}
\label{1.1}
\partial_t \theta  + u \cdot \nabla \theta =0,
\end{equation}
where
\begin{align}
\label{1.2} u=(u_1,u_2) =  \left( \frac{ \partial \psi}{\partial x_2}  , - \frac{\partial \psi}{\partial x_1}
\right),\\
\label{1.3} \theta(x,t) = \Lambda^{2-\alpha} \psi(x,t),
\end{align}
and $\Lambda = \sqrt {-\Delta}$ .
 We denote this system as $active \
scalar \ flow$, where the scalar field $\theta(x,t)$ plays the role
of a temperature field.
In this paper we consider the
evolution of a positive compactly supported initial temperature $\theta_0\in L^{\infty}(\mathbb{R}^2)$
according to \eqref{1.1} for $\alpha \in (0,1)$ . Here we assume that 
the support of $\theta_0$ is contained in the ball centred at the origin of bounded diameter $2R_0$. 
The main aim of this paper is to give an estimate on the growth rate of the diameter of the support of a hot patch of temperature evolving
according to \eqref{1.1}. From the physical point of view the problem 
here considered is the confinement in time of a hot patch of temperature 
in absence of an external forcing, under the quasi-geostrophic approximation.\\
We briefly recall the main results obtained about the growth rate
of the diameter of the support of a vortex patch
in the case of the Euler equation in 2D (that is \eqref{1.3} for 
$\alpha= 0$).
The first result, obtained by Marchioro in \cite{Marb}, gives a bound of $O(t^{1/3})$
on the diameter of the support. An improvement of this result was given by
Iftimie et al. in \cite{Gam}, where a bound of $O\left((t\ln t)^{1/4}\right)$
was proved. For the sake of completness we also recall that a similar estimate was
independently obtained at the same time by Serfati in \cite{Serfati}.\\
In the paper by Iftimie et al. \cite{Gam}, the main result was obtained in two different ways.
Indeed a first proof of Proposition 2.1 in \cite{Gam} was based on the application of sharp estimates on
the mass of vorticity near the boundary by using a mollifier (similarly to the proof given in \cite{Marb}).
The second proof in the Appendix of the same paper, was given by Gamblin and is based on estimates
of higher order momenta
\begin{equation}
m_n(t)=\int_{\mathbb{R}^2}|x|^{4n}\theta(x,t) dx.
\end{equation}
 In this paper we will generalize the last result,
for the estimate of the growth of the support of a temperature patch according to
the $\alpha$-patch model. We will show that the proof is based on non-trivial adaptation of the 
tecniques used in \cite{Gam}, due to the different Green function of the problem.

\section{The physical motivation: confinement of an hot patch of temperature in the SQG model}

Here we give an outline of the physical problem considered in this paper, firstly recalling the derivation of the model equations and their 
meaning.\\
The behavior of the atmosphere in the mid-latitudes is described by the equations of rotating fluids and in many cases 
geophysicists adopt the so-called quasi-geostrophic approximation, that is a perturbative expansion around a uniformly rotating state
(see \cite{Ped}). Under this approximation, the long-scale dynamics of the atmosphere in mid-latitudes is governed
by the balance between the Coriolis force and the pressure gradient, i.e. 
\begin{equation}\label{stre}
\begin{cases}
\displaystyle{fu_2 = -\frac{\partial p}{\partial x_1}} , \\
\displaystyle{f u_1 = \frac{\partial p}{\partial x_2}},
\end{cases}
\end{equation}
where we denoted with $u\equiv (u_1, u_2)$ the velocity field, $p(x,y,t)$ the pressure field and $f$ is the Coriolis parameter (assumed constant).
We observe that, by means of \eqref{stre}, we can consider the pressure field as the stream function in this kind of problems in fluid mechanics.
In order to obtain the evolution equations for the pressure field, a perturbative expansion with respect to the zeroth-order 
geostrophic field should be done. At the first order of this expansion, the equation of conservation of the potential vorticity arise, that is
\begin{equation}
\left(\frac{\partial}{\partial t}+\left(u\cdot \nabla\right)\right)q= 0,
\end{equation}
where 
\begin{equation}\label{PV}
q= \frac{\partial^2\psi}{\partial x_1^2}+ \frac{\partial^2\psi}{\partial x_2^2}+\frac{\partial}{\partial x_3}\left(\left(\frac{f}{N}\right)^2
\frac{\partial\psi}{\partial x_3}\right).
\end{equation}
In equation \eqref{PV}, $q(\cdot)$ is the potential vorticity and $N$ the buoyancy frequency. The further assumption in the mathematical derivation
of the surface quasi-geostrophic equation (SQG) is to take $(f/N)^2= 1$ and
$q= const.$ (for simplicity null), therefore obtaining from \eqref{PV}
\begin{equation}\label{PV1}
\left(\frac{\partial^2}{\partial x_1^2}+
\frac{\partial^2}{\partial x_2^2}+\frac{\partial^2}{\partial x_3^2}\right) \psi = 0.
\end{equation}
Is moreover usual in geophysics to apply the Boussinesq approximation, for
which the vertical derivative of the pressure field is proportional to 
the temperature field $\theta$, that is
\begin{equation}\label{PV2}
\theta\propto \frac{\partial p}{\partial x_3}  = \frac{\partial \psi}{\partial x_3} . 
\end{equation}
Moreover, from the energy balance law, we have
that also the temperature field is advected by the fluid
\begin{equation}
\left(\partial_t+\left(u\cdot \nabla\right)\right)\theta=0.
\end{equation}
By using equations \eqref{PV1} and \eqref{PV2} we finally have
\begin{equation}\label{PV3}
\theta= \frac{\partial}{\partial x_3}\psi= \left(-\Delta\right)^{1/2}\psi,
\end{equation}
where $$\Delta = \frac{\partial^2}{\partial x_1^2}+\frac{\partial^2}{\partial x_2^2}.$$
We recall that the operator appearing in \eqref{PV3} is a pseudo-differential operator, defined by the Fourier symbol as
\begin{equation}
(-\hat{\Delta})^{1/2}\psi(k) = |k| \hat{\psi}(k)
\end{equation}
and it is the operator mapping Neumann to Dirichlet boundary value problems (for a discussion about this point see e.g. \cite{caffarelli}).\\
Assuming that the fluid is incompressible, i.e. that the velocity field is divergence free, we finally obtain the
basic equations of the SQG model
\begin{equation}
\begin{cases}
\partial_t \theta +(u\cdot \nabla)\theta = 0 ,\\
\nabla \cdot u= 0,
\end{cases}
\end{equation}
where $\theta = (-\Delta)^{1/2}\psi$. \\
Starting from the seminal paper of Constantin, Majda and Tabak \cite{Const1} about the singular front formation 
in a model for a quasi-geostrophic flow, a great number of works about the SQG equation and its generalizations 
with dissipation appeared in literature (see e.g. \cite{Chae},\cite{Inv4}, \cite{Li},\cite{Const2},\cite{Rodrigo}
and the references therein). 
The modified SQG equation was recently studied in different papers as 
a family of active scalar flows interpolating the 2D Euler and SQG equations.
Indeed, we observe that, equations \eqref{1.1}-\eqref{1.2}-\eqref{1.3}, for $\alpha =0$ describe the
flow of an ideal incompressible fluid according to the Euler equation in two
dimensions in which $\theta$ has the meaning of vorticity, while for
$\alpha =1$ the Surface Quasi-Geostrophic equation (SQG) is recovered in which
$\theta$ has the meaning of temperature. This last system has
relevance in geophysics as just shown (see for instance \cite{Ped}). The modified SQG
equation was firstly studied by Pierrhaumbert et.al
\cite{Held}, in relation to the spectrum of turbulence
in 2-d, beginning the research on $\alpha-turbulence$ (see also \cite{Held1}). On the other
hand, from a mathematical point of view, this dynamical system was
deeply studied (see \cite{Const1},\cite{Const2}) in relation to the
problem of formation of singular fronts. The authors observed a
formal analogy with the 3D Euler  equation and studied the singular
behavior of the solutions in the more suitable framework of the 2D
quasi-geostrophic equation.\\
The patch problem for SQG-type equations has been considered by different authors in recent papers and
leads to nontrivial problems. Local existence results were proved by Gancedo in \cite{Gancedo}
 for the SQG and modified SQG patches with boundaries which are $H^3$ closed curves.
Computational studies have been also performed in \cite{Cordoba1} and \cite{Mancho}, suggesting
a finite time singularity with two patches touching each other. 
Finally, in the recent paper \cite{kisel}, Kiselev et al. have studied
the modified SQG patch dynamics in domains with boundaries. Some results in this direction have been also proved by Cavallaro et al. in \cite{Cav}
where the dynamics of $N$ concentrated temperature fields has been considered.\\
As far as we know, a detailed analysis of the 
growth in time of the radius of a single patch of temperature evolving according to the modified SQG equation in the whole plane is still 
missing in the literature, this will be the main object of our paper.\\
We also briefly add some considerations on the modified SQG that
will be useful in the next sections.\\
It is possible to introduce a weak form of eq.s (1.1)-(1.4):
\begin{equation}
\label{2.1} \frac{d}{dt}\theta[f] = \theta[u \cdot \nabla f]+
\theta[\partial_t f],
\end{equation}
where  $f(x,t)$ is a bounded smooth function and
\begin{equation}
\label{2.2} \theta [f] = \int dx \ \theta(x,t) \ f(x,t) \ .
\end{equation}
In \eqref{2.1} the velocity field $u$ is given by
\begin{equation}\label{u}
u(x,t)=\int K(x-y)\theta(y,t)dy,
\end{equation}
where $K(x-y)= \nabla^\perp G(x-y)$, being $\nabla^\perp =
(\partial_{x_2},-\partial_{x_1})$ and $G(x)$ the Green function of
$(-\Delta)^{(1-\alpha/2)}$ with vanishing boundary condition at
infinity.

Moreover the divergence-free of the velocity field and  the fact
that $\theta$ is transported by the flow show that the Lebesgue
measure and the maximum of $\theta$ are conserved during the motion.

\bigskip

\noindent We remember that the Green function of
$(-\Delta)^{(1-\alpha/2)}$ in $\R^2$ with vanishing boundary
conditions at infinity is (see for instance \cite{Green}):
\begin{equation}
\label{2.8} G(r)= \psi (\alpha) r^{-\alpha} \ \ , \ \ \psi(\alpha) =
-\frac{1}{\pi 2^{(2-\alpha)}}
\frac{\Gamma(\alpha/2)}{\Gamma(\frac{2-\alpha}{2})} \ \ , \ \
r=\sqrt{x^2_1+x^2_2} \ ,
\end{equation}
where $\Gamma(\cdot)$ denotes the Euler Gamma function.\\
Finally, we recall, for the utility of the reader, that several
quantities are conserved in the time evolution for the
$\alpha$-patch model (see e.g. \cite{Cav}), that is
\begin{itemize}
\item the total mass
\begin{equation}
\nonumber\int \theta(x,t) dx= \int \theta_0(x) dx= m_0.
\end{equation}
\item the maximum norm
\begin{equation}
\|\theta(x,t)\|_{L^{\infty}}= \|\theta_0(x)\|_{L^{\infty}} = M_0.
\end{equation}

\item the center of mass
\begin{equation}
\nonumber \frac{1}{m_0}\int x \ \theta(x,t) dx= \frac{1}{m_0}\int
x \ \theta_0(x) dx= c_0.
\end{equation}
\item the moment of inertia
\begin{equation}
\nonumber \int |x|^2 \ \theta(x,t) dx= \int |x|^2 \ \theta_0(x) dx= i_0.
\end{equation}
\end{itemize}
These conservation laws play a central role for the proof of the
main theorem in the next section.

The aim of this paper is to rigorously study the confinement of an hot patch of temperature in absence of an external forcing .
We are therefore interested in proving some rigorous results on the concentration of a temperature field under the 
quasi-geostrophic approximation, in the more \textit{flexible} formulation of the modified SQG equation in which the parameter
$\alpha$ (that appears in the fractional Laplacian of equation \eqref{1.3}) describe the difference in the behavior
of the air mass between the description given by the ideal incompressible fluid (according to the Euler equation)
 and the one provided by the SQG equation.
In particular in the next section we generalize some rigorous results previously proved by Iftimie et al. in \cite{Gam}
for the Euler equation in the plane, in order to estimate the growth in time of a hot patch of temperature
according to the $\alpha$-patch model.

\section{Estimates on the growth of the support of a temperature patch}

\label{sec:3} We now give an estimate on the growth of the
support of a temperature patch, assuming that the support of the initial temperature field $\theta_0$ is contained
in the ball centered at the origin of radius $R_0$.
\begin{theorem}
Let $\theta(x,t)$ be a solution of the active scalar flow equation with a positive compactly supported initial temperature
$\theta_0\in L^{\infty}(\mathbb{R}^2)$. There exists a constant $C_0\equiv C_0(i_0, R_0, m_0, M_0)$ such that, for every $t\geq 0$ the support of 
$\theta(x,t)$ is contained in the ball $|x|< 4R_0+C_0\left[t\ln(2+t)\right]^{\frac{1}{4+\alpha}},$ for any $\alpha \in(0,1)$.
\end{theorem}

\smallskip

\textit{Proof}\\
The proof of the theorem is strongly based on the use of the conserved quantities and in particular of 
the moment of inertia $i_0$ and the total mass $m_0$. The main difference and technical difficulty with respect to the case $\alpha =0$ is due to the different Green function given by \eqref{2.8}. \\
Without loss of generality, we will take the center of mass located at the origin. We should prove that the radial component of the velocity field 
satisfies the following estimate
\begin{equation}\label{imp}
\bigg|\frac{x}{|x|}\cdot u(x,t)\bigg|\leq \frac{C_0}{x^{3+\alpha}},
\end{equation}
for all $|x|\geq 4d_0+C_0\left[t\ln(2+t)\right]^{\frac{1}{4+\alpha}}$.\\
The radial component of the velocity field in the case of the active scalar 
flow is given by
\begin{equation}
\frac{x}{|x|}\cdot u= \frac{x}{|x|}\cdot\int\nabla^{\perp}G(x,y)\theta(y)dy= 
\int\frac{x}{|x|}\cdot \frac{(x-y)^{\perp}}{|x-y|^{2+\alpha}}\theta(y)dy,
\end{equation}
where $(x_1,x_2)^{\perp}=(-x_2,x_1)$.
Let us divide the last integral in two regions:
\begin{align}
\nonumber &A_1:= |x-y|<\frac{|x|}{2}\\
\nonumber &A_2:= |x-y|\geq\frac{|x|}{2}.
\end{align}
The integral over the region $A_1$ is simply bounded by 
\begin{equation}
C\int_{A_1}\frac{\theta(y)}{|x-y|^{\alpha+1}}dy.
\end{equation}
By using the fact that $x\cdot (x-y)^{\perp}= -x\cdot y^{\perp}$ and the 
conservation of the center of the mass (that is located at the origin), we have
\begin{align}
\nonumber&\int_{A_2}\frac{x}{|x|}\cdot \frac{(x-y)^{\perp}}{|x-y|^{\alpha+2}}\theta(y)
dy=- \int_{A_2}\frac{x\cdot y^{\perp}}{|x||x-y|^{2+\alpha}}\theta(y)dy\\
&=-\int_{A_2}\frac{x\cdot y^{\perp}}{|x|}\left(\frac{1}{|x-y|^{2+\alpha}}-\frac{1}{|x|^{2+\alpha}}\right)\theta(y)dy+ \int_{A_1}\frac{x\cdot y^{\perp}}{|x|^{3+\alpha}}\theta(y)dy
\end{align}
The first integral can be bounded as follows
\begin{equation}
\bigg|\int_{A_2}\frac{x\cdot y^{\perp}}{|x|}\left(\frac{1}{|x-y|^{\alpha+2}}-\frac{1}{|x|^{\alpha+2}}\right)\theta(y) \ dy\bigg|\leq \frac{const \ i_0}{|x|^{3+\alpha}},
\end{equation}
where we have used the conservation of the moment of inertia.\\
Regarding the integral over the region $A_1$, being $|x-y|<\frac{|x|}{2}$,
then we have that $|y| \leq \frac{3|x|}{2}$ and therefore
\begin{equation}
\bigg|\int_{A_1}\frac{x\cdot y^{\perp}}{|x|^{3+\alpha}} \ \theta(y)dy\bigg|\leq const. \ \int_{A_1}\frac{\theta(y)}{|x-y|^{\alpha+1}}dy.
\end{equation}
At this stage, we have obtained the following estimate for the radial component of the velocity field
\begin{equation}
\bigg|\frac{x}{|x|}\cdot u(x)\bigg|\leq \frac{const. \ i_0}{|x|^{3+\alpha}}+const. \ \int_{A_1}\frac{\theta(y)}{|x-y|^{\alpha+1}}dy
\end{equation}
We now introduce the following Lemma
\begin{lemma}\label{lem}
Let $S\subset \mathbb{R}^2$ and $h: S\rightarrow \mathbb{R}^+$ a function belonging to $L^1(S)\bigcap L^{\infty}(S)$. Then
\begin{equation}
\int_S\frac{h(y)}{|x-y|^{\beta}}dy\leq C \ \|h(y)\|_{L_1}^{1-\frac{\beta}{2}}\|h(y)\|_{L^{\infty}}^{\beta/2},
\end{equation} 
with $\beta\in (0,2)$.
\end{lemma}
The proof is essentially the same of Lemma 2.1 in \cite{Gam}, details of the proof are given in the Appendix.\\
By using the fact that 
\begin{equation*}
\{y:|x-y|<|x|/2\}\subset\{y: |y|>|x|/2\}
\end{equation*}
and the Lemma \ref{lem}  we have the following estimate on 
the radial component of the velocity
\begin{equation}
\bigg|\frac{x}{|x|}\cdot u\bigg|\leq \frac{C}{|x|^{3+\alpha}}+
C \ M_0^{\frac{\alpha+1}{2}}\left(\int_{|y|>|x|/2}\theta(y)dy\right)^{\frac{1-\alpha}{2}}.
\end{equation}
In the next proposition we prove that also the last integral is of order 
$O(|x|^{-3-\alpha})$ for $|x|$ large, so that inequality \eqref{imp} holds and the proof of the theorem is achieved.

\smallskip

\begin{proposition}
There exists a constant $C_0\equiv C_0(i_0,R_0, m_0,M_0)$ such that, for any $k>0$
\begin{equation}
\int_{|y|>|x|/2}\theta(y,t)dy\leq \frac{C_0}{|x|^{k}},
\end{equation}
for all $|x|>4R_0 +C_0\left(t \ \ln(2+t)\right)^{\frac{1}{4+\alpha}}$.
\end{proposition}

\textit{Proof.}\\
The proof of this proposition is based on the following technical Lemma:

\begin{lemma}
Let be
\begin{equation}
m_{n, \alpha}(t)=\int |x|^{(4+\alpha)n}\theta(x,t)dx,
\end{equation}
then, there exists a constant $C_0$ such that, for any $n\geq 1$
\begin{equation}
m_{n, \alpha}(t)\leq m_0\left(R_0^{4+\alpha}+C_0 i_0n^{1+\alpha}t\right)^n.
\end{equation}
\end{lemma}
We postpone the proof of this Lemma and use it to prove Proposition 3.1.

Let us suppose that 
\begin{equation}
r^{4+\alpha}\geq 2[R_0^{4+\alpha}+C_0 i_0k t\ln(2+t)],
\end{equation}
taking $n$ such that 
\begin{equation}
k\ln(2+t)-1< n^{\alpha+1}\leq k\ln(2+t),
\end{equation}
then we have from Lemma 4.2 that 
\begin{align}
\nonumber \int_{|x|\geq r}\theta(x,t)dx\leq \frac{m_{n,\alpha}(t)}{r^{4(n+\alpha)}}
\leq \frac{m_0}{r^{k}}\frac{\left(R_0^{4+\alpha}+C_0 i_0n^{\alpha+1}t\right)^n}{r^{4n-k+n\alpha}}\leq const. \ \frac{m_0}{r^k},
\end{align}
when $r^{4+\alpha}\geq 2[R_0^{4+\alpha}+C_0 \ i_0k t \ln(2+t)]$, as claimed.

\bigskip

\textit{Proof of Lemma 4.2.}\\
By using \eqref{1.1} and \eqref{u} and integrating by parts, we have that
\begin{align}
m_{n, \alpha}'(t)&=\int |x|^{(4+\alpha)n}\partial_t \theta(x,t)dx
=\int |x|^{(4+\alpha)n}u_r\cdot\nabla \theta(x,t)dx\\
\nonumber & =\frac{(4+\alpha)n}{2\pi}\int\int \frac{x\cdot(x-y)^{\perp}}{|x-y|^{2+\alpha}}
|x|^{4n+\alpha n-2}\theta(x,t)\theta(y,t)dxdy.
\end{align}
We define
\begin{equation}
\mathcal{K}(x,y)=\frac{x\cdot(x-y)^{\perp}}{|x-y|^{2+\alpha}},
\end{equation}
so that, we can write in a compact way
\begin{equation}
m_{n, \alpha}'(t)= \frac{(4+\alpha)n}{2\pi}\int\int \mathcal{K}(x,y)|x|^{4n+\alpha n-2}\theta(x,t)\theta(y,t)
dxdy.
\end{equation}
We consider the following partition
\begin{align}
\nonumber &A_1=\{(x,y):|y|\leq\left(1-\frac{1}{2n}\right)|x|\}\\
\nonumber &A_2=\{(x,y):\left(1-\frac{1}{2n}\right)|x|\leq |y|\leq \left(1-\frac{1}{2n}\right)^{-1}|x|\}\\
\nonumber &A_3=\{(x,y):|x|\leq \left(1-\frac{1}{2n}\right)|y|\}.
\end{align}
Then, we have that
\begin{equation}
\frac{d}{dt}m_{n, \alpha}(t)= \sum_{i=1}^3 \frac{(4+\alpha)n}{2\pi}\int\int_{A_i} \mathcal{K}(x,y)
|x|^{4n+\alpha n-2}\theta(x,t)\theta(y,t)
dxdy.
\end{equation}
We first consider the region $A_1$, where $|x-y|\geq |x|/2n$
and we have the following
bound
\begin{equation}
|\mathcal{K}(x,y)|\leq \frac{|x|}{|x-y|^{\alpha+1}} \leq const.\; \frac{n^{\alpha+1}}{|x|^{\alpha}},
\end{equation}
such that

\begin{align}
\nonumber &\frac{(4+\alpha)n}{2\pi}\int\int_{A_1} \mathcal{K}(x,y)
|x|^{4n+\alpha n-2}\theta(x,t)\theta(y,t)
dxdy\leq const.\; n^{2+\alpha}\int\int_{A_1}
|x|^{(4+\alpha)n-\alpha-2}\theta(x,t)\theta(y,t) dxdy \\
& \leq const.\; \frac{n^{4+\alpha}}{(1-2n)^2}\int\int_{A_1}
|x|^{(4+\alpha)(n-1)}|y|^2\theta(x,t)\theta(y,t) dxdy \leq const.\; n^{2+\alpha} i_0 m_{n-1,\alpha}(t),
\end{align}
where we used also the fact that $|y|^2/|x|^2\geq (\frac{1}{2n}-1)^2$ in $A_1$ and the conservation
of the moment of inertia.\\
In the set $A_3$, we have that $|x-y|\geq|y|/2n$ and we can
obtain the following estimate
\begin{equation}
|\mathcal{K}(x,y)|\leq \frac{|x|}{|x-y|^{\alpha+1}}
\leq const. \; n^{\alpha+1}\frac{\left(1-\frac{1}{2n}\right)}{|y|^{\alpha}}\leq
 const.\; n^{\alpha+1}\frac{|y|^2}{|x|^{\alpha+2}},
\end{equation}
where we also used the fact that in $A_3$, $(1-\frac{1}{2n})\frac{|y|}{|x|}\geq 1$.
It follows that, also in this case,
\begin{align}
\frac{(4+\alpha)n}{2\pi}\int\int_{A_2} \mathcal{K}(x,y)
|x|^{4n+\alpha n-2}\theta(x,t)\theta(y,t)
dxdy \leq const. n^{2+\alpha} i_0 m_{n-1,\alpha}(t).
\end{align}
Finally, the contribution related to the set $A_2$
 can be bounded in a similar way by
a term proportional to $m_{n-1, \alpha}(t)$. Thus we have the following estimate
\begin{equation}
\frac{d}{dt}m_{n, \alpha}(t)\leq C_0 i_0 n^{2+\alpha} m_{n-1,\alpha}(t).
\end{equation}
We now rearrange the last term
\begin{equation}
\int |x|^{(4+\alpha)(n-1)}\theta(x,t)dx=\int |x|^{(4+\alpha)(n-1)}\theta^{1/n}(x,t)
\theta^{1-1/n}(x,t)dx,
\end{equation}
and applying the Holder's inequality, we finally obtain
\begin{equation}
\frac{d}{dt}m_{n,\alpha}(t)\leq C_0i_0n^{2+\alpha} m_0^{1/n}m_{n,\alpha}^{1-1/n}(t),
\end{equation}
and integrating we obtain the claimed result.

\begin{remark}
As a final remark, we observe that uniqueness result for the global solution of the modified SQG equation with an initially positive compactly supported datum is still missing and therefore we are forced to state our main result in a weak form, that is a rigorous estimate of the growth of the support of a solution of the governing equation.
\end{remark}

\smallskip

\textbf{Acknowledgments}\\
We wish to thank Carlo Marchioro (Rome) for fruitful discussions.

\section{Appendix: Proof of Lemma 3.2}
For the sake of completness, we here provide the detailed proof of Lemma 3.2. In particular, we here recall a more general Lemma (see \cite{dra})
that includes Lemma 3.2 as a special case
\begin{lemma}
Let $\beta \in (0,2)$, $S\subset \mathbb{R}^2$ and $h: S\rightarrow \mathbb{R}^+$ a function belonging to $L^1(S)\bigcap L^{p}(S)$, with 
$p>\frac{2}{2-\beta}$. Then
\begin{equation}
\int_S\frac{h(y)}{|x-y|^{\beta}}dy\leq C \ \|h(y)\|_{L_1(S)}^{\frac{2-\beta-2/p}{2-2/p}}\|h(y)\|_{L^{p}(S)}^{\frac{\beta}{2-2/p}},
\end{equation} 
with $\beta>0$.
\end{lemma}

\smallskip

\textit{Proof}\\

Let $k>0$, by using the H\"older inequality, we have that
\begin{align}
\nonumber &\int_S \frac{h(y)}{|x-y|^\beta}= \int_{S\bigcap \{|x-y|>k\}}\frac{h(y)}{|x-y|^\beta}dy+ \int_{S\bigcap \{|x-y|<k\}}\frac{h(y)}{|x-y|^\beta}dy\\
\nonumber & \leq \frac{\|h\|_{L^1(S)}}{k^\beta}+
\|h\|_{L^p(S)}\bigg\|\frac{1}{|x|^\beta}\bigg\|_{L^{\frac{p}{p-1}}(|x|\leq k)}\\
\nonumber & = \frac{\|h\|_{L^1(S)}}{k^\beta}+C \|h\|_{L^p(S)} k^{2-\beta-2/p}.
\end{align}
and by taking $k = \left(\|h\|_{L^1(S)}\|h\|_{L^p(S)}^{-1}\right)^{\frac{1}{2-2/p}}$, 
we conclude the proof.\\
Finally observe that the Lemma 3.2 is a particular case, for $h \in L^1(S)\bigcap L^{\infty}(S)$.

\end{document}